
\magnification=1200
\baselineskip=13.75pt
\vsize=9.55truein
\hsize=6.25truein
%
%

\font\bigbf=cmbx10 at 14pt

\font\smallbf=cmbx10 at 9pt
\font\smallrm=cmr10 at 9pt

%
%
%
%
\def\qed{${\vcenter{\vbox{\hrule height .4pt
           \hbox{\vrule width .4pt height 4pt
            \kern 4pt \vrule width .4pt}
             \hrule height .4pt}}}$}
\def\mqed{{\vcenter{\vbox{\hrule height .4pt
           \hbox{\vrule width .4pt height 4pt
            \kern 4pt \vrule width .4pt}
             \hrule height .4pt}}}}
%
%
\font\bbb=msbm10 at 10pt  

\def\HH{\hbox{\bbb H}}

\def\RR{\hbox{\bbb R}}
\def\CC{\hbox{\bbb C}}
\def\ZZ{\hbox{\bbb Z}}

\def\QQ{\hbox{\bbb Q}}

\def\GG{\hbox{\bbb G}}

\newcount\refCount
\def\newref#1 {\advance\refCount by 1
\expandafter\edef\csname#1\endcsname{\the\refCount}}
\newref ABRAMOVICH
\newref ANDREBOOK
\newref ANDRE
\newref AX
\newref BILUPARENT
\newref HEIGHTS
\newref CGMM
\newref FALTINGS
\newref FALTINGSBIG
\newref FREY
\newref HPA
\newref HPB
\newref KLINGLERYAFAEV
\newref KUHNE
\newref LANG
\newref LANGALGEBRA
\newref LAURENT
\newref LENSTRA
\newref LIARDET
\newref MANN
\newref MASSER
\newref MEREL
\newref NAJMAN
\newref OORT
\newref PELLARIN
\newref PILAOAO
\newref PILAEMS
\newref PILAICM
\newref PTAG
\newref PTAXJ
\newref PILAWILKIE
\newref PILAZANNIER
\newref PINK
\newref POONENGONALITY
\newref RAYNAUDA
\newref RAYNAUDB
\newref REBOLLEDO
\newref RISMAN
\newref SARNAKADAMS
\newref SIEGEL
\newref TSIMERMAN
\newref ULLMOYAFAEV
\newref WILES
\newref WILKIE
\newref ZANNIER
\newref ZILBER
\newref ZOGRAF

\newcount\derandomCount
\def\derandom#1 {\advance\derandomCount by 1
\expandafter\edef\csname#1\endcsname{\the\derandomCount}}
\rightline{20150820}

\bigskip
\bigskip

\centerline{\bigbf On a modular Fermat equation}

\bigskip
\centerline{J. Pila}


\bigskip
\bigskip
\noindent
{\narrower 
\baselineskip=11pt
{\smallbf Abstract.\/} {\smallrm We consider some
diophantine problems suggested by the analogy between 
multiplicative groups and powers of 
the modular curve in problems of ``unlikely intersections''.
We prove a special case of the 
Zilber-Pink conjecture for curves.}

\medskip
\noindent
{\smallbf 2010 Mathematics Subject Classification:\/} 
{\smallrm 11G18, 11D41, 03C64 (primary).\/}

}

\bigskip
\bigskip
\bigskip
\centerline{\bf 1. Introduction}
\medskip


To motivate the title problem, we recall some classical
diophantine statements.
We identify (algebraic) varieties with their sets of 
complex points. Thus, in particular,  
$\GG_{\rm m}=\GG_{\rm m}(\CC)=\CC^\times$
is the multiplicative group of non-zero complex numbers, and
$Y(1)=Y(1)(\CC)=\CC$ is the moduli space parameterising 
elliptic curves defined over $\CC$,
up to isomorphism over $\CC$, by their $j$-invariant.

The  Multiplicative Manin-Mumford
conjecture (MMM; a theorem Laurent
[\LAURENT], see also [\MANN, \SARNAKADAMS])
concerns the distribution of {\it torsion points\/} in a subvariety
$V\subset \GG_{\rm m}^k$. These are the torsion points in the
group, namely the points of the form $(\zeta_1,\ldots, \zeta_k)$
where $\zeta_i\in \CC^\times$ are roots of unity.
MMM states that the torsion points contained in $V$ are 
contained in a finite number of {\it torsion cosets\/} contained in $V$.
Torsion cosets are the cosets of subtori by torsion points; 
otherwise expressed, they are the 
irreducible components of subvarieties defined by systems of
multiplicative relations, that is relations of the form 
$x_1^{a_1}\ldots x_k^{a_k}=1$. Thus a torsion point is precisely
a torsion coset of dimension zero.
The ``original'' Manin-Mumford conjecture (MM) is the same statement
for a subvariety of an abelian variety, in which torsion cosets
are cosets of abelian subvarieties by torsion points.
MM is a theorem of Raynaud
[\RAYNAUDA,\RAYNAUDB].

The Andr\'e-Oort conjecture [\ANDREBOOK, \OORT] 
was partly motivated by an informal analogy with MM.
It concerns the distribution of {\it special points\/} in a subvariety
$V$ of a {\it Shimura variety\/} $X$,
and is now ``almost'' fully proved [\KLINGLERYAFAEV, 
\ULLMOYAFAEV, \TSIMERMAN].
For cartesian powers of the modular curve 
it is a theorem (``Modular Andr\'e-Oort''; MAO)
proved in [\ANDRE, \PILAOAO]. MAO states that,
for $V\subset Y(1)^k$, the {\it special points\/} of $Y(1)^k$ 
contained in $V$ are contained in finitely many 
{\it special subvarieties\/} of $Y(1)^k$ 
contained in $V$. 
The special subvarieties of $Y(1)^k$ 
are the irreducible components of subvarieties defined 
by systems of modular relations,
that is relations of the form $\Phi_{N_{ij}}(x_i,x_j)=0$, 
where $\Phi_N$ are the classical modular polynomials.
A special point is precisely a special subvariety of dimension zero.
See \S6 for a more careful definition for $Y(1)^k$, 
and \S2 for $Y(1)^2$.

The two statements are unified within Pink's version [\PINK]
of what is now known as the Zilber-Pink conjecture (ZP).
See also [\ZILBER, \ZANNIER, \PILAICM] 
for the general formulation of this far-reaching conjecture, 
which is very much open, and
\S7 for the statement in $Y(1)^k$.
ZP governs the interaction between a subvariety
$V$ of a mixed Shimura variety $X$, and the collection of 
{\it special  subvarieties\/} of $X$ (see [\PINK]). 
In $\GG_{\rm m}^k$, the special subvarieties
are the aforementioned torsion cosets.
Thus, within ZP, MMM and MAO are analogues in a strict sense,
and modular relations are analogues of multiplicative ones.

In the multiplicative setting,  the Multiplicative Mordell-Lang 
conjecture (MML;  a theorem of Laurent [\LAURENT])
generalises MMM.
Let us state it in the special case of the variety 
$V\subset \GG_{\rm m}^2$ defined by $u+v=1$
(the {\it unit equation\/}): there are only finitely 
many solutions to $u+v=1$ when $u,v$ are restricted to the 
division group of a finitely 
generated subgroup of $\CC^\times$. 
Important special cases, for finitely generated subgroups
of algebraic, or even rational, numbers were 
established in fundamental work of Siegel, Mahler, 
Lang, and Liardet; see [\SIEGEL, \LANG, \LIARDET, \HEIGHTS].

The modular analogue of this
statement  (``Modular Mordell-Lang'')
is proved, in general form, in [\HPA, \PILAEMS]. In the special
case it asserts that there are only finitely many solutions to $u+v=1$
when $u,v$ are restricted to finitely many {\it Hecke orbits\/}
(or are special points). The Hecke orbit of
$x\in \CC$ is $\{y\in \CC: \exists N\  \Phi_N(x,y)=0\}$.
It is the set of $j$-invariants of elliptic curves which are {\it isogenous\/}
to the one with $j$-invariant $x$.

Now we observe that Fermat's Last Theorem 
(FLT; theorem of Wiles [\WILES])
may also be expressed in these terms: 
it asserts that $u+v=1$ has no solutions
for $u,v\in \QQ^{\times n}$ when $n\ge 3$. It seems not to have been 
observed that the condition on $u,v$ fits naturally into the
multiplicative group setting: they are required to be in the subgroup 
consisting of $n$th powers of rational numbers.  
The modular analogue of $u=x^n$ is 
$\Phi_n(x,u)=0$. Generalising a little, we are led to investigate the rational
solutions $x,y$ of the system
$$
\Phi_N(x,u)=0, \quad \Phi_M(y,v)=0,\quad u+v=1,\quad N,M\ge 1.
\eqno{(*)}
$$

This is the ``modular Fermat equation'' of the title.
Given $N, M$ one may eliminate $u,v$ in $(*)$ to find that $(x,y)$
lie on some algebraic curve (possibly reducible) $V_{N,M}$.
The strict analogue of FLT would take $N=M$, but this plays no role
for us.  We prove the following partial analogue of ``asymptotic''
FLT. It asserts that there are no rational points
on any of the curves $V_{N,M}$ with large prime $\max\{N,M\}$.
We have nothing to say about possible solutions for
small $N, M$.

\medskip
\noindent
{\bf 1.1. Theorem.\/} {\it  There exists $L$ such that 
$(*)$ 
has no solutions with $x,y\in\QQ$ for which
$\max\{N,M\}\ge L$ and  $\max\{N,M\}$ is a prime number.
\/}
\medskip

Our proof of this theorem uses a variant of the o-minimality and 
point-counting strategy which has been used over recent years 
to prove  various cases of the Andr\'e-Oort (and Zilber-Pink) 
conjecture, using the Counting Theorem of Pila-Wilkie [\PILAWILKIE].
The strategy was  originally proposed by Zannier in the context of the
Manin-Mumford conjecture (see [\PILAZANNIER]), where it relies
on torsion points having high degree (relative to their order).
For Andr\'e-Oort, the strategy
depends on special points having high degree over $\QQ$
in a suitable sense (see [\PTAG, \TSIMERMAN]).
Our results here likewise depend on $\QQ(u,v)$ having 
large degree over $\QQ$ (relative to $\max\{N,M\}$, in a 
sense made precise below). 
The applicability of the Counting Theorem in these settings
relies ultimately on the result of Wilkie [\WILKIE] that the real
exponential function gives rise to an {\it o-minimal structure\/}.
Before going further 
into the specifics, let us observe that this method has no purchase
for FLT or Mordell-Lang type problems, simply because when
$u, v$ are in a group generated by rational numbers, or a finitely 
generated group, $[\QQ(u,v):\QQ]$ is bounded.

We can remove the primality condition on $\max\{N,M\}$ conditionally 
on a special case 
of a statement (``GO1'', see \S8) formulated in Habegger-Pila [\HPB]. 
Consider $x,y\in \overline{\QQ}$ such that the elliptic curves 
$E_x, E_y$,
whose $j$-invariants are $x,y$, are related by a cyclic isogeny
of degree $N$. So $\Phi_N(x,y)=0$ for the modular 
polynomial $\Phi_N$.
If $x,y$ are not special (see \S2) 
then $N$ is unique. 

\medskip
\noindent
{\bf 1.2. Strong Galois-Orbit Hypothesis (SGH).\/} There exist  
$c,\delta>0$ such that if $(x,y)\in {\overline{\QQ}}^2$ are not
special and $\Phi_N(x,y)=0$, then
$$
[\QQ(x,y):\QQ]\ge c N^\delta.
$$

\medskip

The plausibility of this conjecture is discussed briefly in \S3. 
Essentially, it is on a par with expectations for the
best dependence in the
Strong Uniform Boundedness Conjecture (Merel's theorem [\MEREL]).
In \S8 we show that SGH is the essential case of the statement 
GO1 alluded to above, and that, in view of [\HPB, \PTAXJ], 
it implies the full Zilber-Pink
conjecture for $Y(1)^k$ (see \S7 for the statement).
We need just the special case of 1.2 for $x\in\QQ$
to prove an unrestricted version of 1.1.

\medskip
\noindent
{\bf 1.3. Theorem.\/} {\it Assume SGH for $x\in\QQ$. Then 
there exists $L$ such that 
$(*)$ 
has no solutions with $x,y\in\QQ$ for $\max\{N,M\}\ge L$.
\/}

\medskip

The reason we are able to prove Theorem 1.1 is that SGH for
$x\in \QQ$ and $N$ a prime number follows from recent results
of Najman [\NAJMAN]. His results are more precise, but imply
in particular that if $\Phi_N(x,y)=0$ with $x\in\QQ$ and $N\ge 41$ 
a prime then 
$$[\QQ(y):\QQ]\ge N/3.
$$

Though very much in the spirit of ``unlikely intersections'', the 
conclusion of 1.3 is seemingly not a consequence of the 
Zilber-Pink conjecture, because rational points in $Y(1)^2$
are neither special
nor contained in finitely many Hecke orbits. Likewise,
FLT is not a consequence of MML because $\QQ^{\times n}$
is not finitely generated.

In \S\S4, 5, 6 we consider generalisations. We can prove
analogues of 1.1 and 1.3 for more general curves
and higher-dimensional varieties in $Y(1)^k$.
These suggest the formulation of analogous conjectures in 
the multiplicative setting which generalise (asymptotic) FLT. 
Our methods cannot address them, but we prove (Theorem 6.4) 
the analogue of one of our main conjectures (5.4), 
for the {\it inverse\/} Fermat equation. This would seem to add 
credence to the conjectures since
$u^n=x$ is also an analogue of $\Phi_n(x,u)=0$.
All our conjectures for curves in \S4 are implied 
by the $abc$ conjecture.

In \S7 and \S8 we study the relationship between SGH and
statements formulated in [\HPB]. We observe that,
if $x,y$ are non-algebraic points with $\Phi_N(x,y)=0$, then 
the large gonality of modular curves ([\ZOGRAF, \ABRAMOVICH])
implies that we get a high 
extension degree even over finitely generated fields. 
Note that gonality growth of some positive power of $N$
is necessary if SGH is to be true.
This enables us to prove a special case of the Zilber-Pink 
conjecture for curves, a counterpart to the result
of [\HPA].

\medbreak
\noindent
{\bf 1.4. Theorem.\/} {\it Let $V\subset Y(1)^3$ be a curve which is
not defined over $\overline{\QQ}$. Then the Zilber-Pink
conjecture holds for $V$.\/}

\medbreak

Note that if $V$ as in 1.4 is not contained in any proper subvariety
of $Y(1)^3$ defined over $\overline{\QQ}$ then the conclusion
follows from the main result of 
Chatzidakis-Ghioca-Masser-Maurin [\CGMM]. We will use this 
in extending the above result to curves in $Y(1)^k$ provided that 
no image under a coordinate projection to $Y(1)^3$ is defined over 
$\overline{\QQ}$.

In our proofs, the Galois and gonality results mentioned (which
show that the points in question have ``many'' conjugates),
are opposed to upper bounds for rational points on suitable
sets definable in an o-minimal structure. This basic strategy
has been used in many problems along these lines.
A new feature here is that the proofs use a family of definable sets,
and rely on uniformity in the Counting Theorem.

\bigbreak
\centerline{\bf 2. Proof of Theorems 1.1 and 1.3}

\medskip

A {\it special point\/} in $\CC$, also known as a {\it singular modulus\/}, 
is the $j$-invariant of a CM 
elliptic curve. Equivalently, it is a number $\sigma=j(\tau)$ where
$\tau\in \HH$ is a quadratic point ($[\QQ(\tau):\QQ]=2$). 
Here $\HH$ is the complex upper half-plane and 
$j:\HH\rightarrow \CC$ is the  elliptic modular function.

\medbreak
\noindent
{\bf 2.1. Definition.\/} A {\it special subvariety\/} of $\CC^2$ is one of
the following: $\CC^2$ itself; a modular curve $T_N$ defined by $\Phi_N(x,y)=0$; a line
$x=\sigma$ or $y=\sigma$ where $\sigma$ is a singular modulus; or
a point $(\sigma, \sigma')$ where $\sigma, \sigma'$ are singular moduli
(a {\it special point\/} of $\CC^2$). {\it Weakly special subvarieties\/} 
include the above,
all horizontal and vertical lines, and all points.

\medbreak
\noindent
{\bf 2.2. Proof of Theorem 1.1.\/} Let $F\subset\HH$ be the standard fundamental domain for the action of ${\rm SL}_2(\ZZ)$ on $\HH$
by Mobius transformations. The  restriction $j: F\rightarrow \CC$
of the elliptic modular function is  definable 
in the o-minimal structure $\RR_{\rm an\ exp}$.

Define the following family of sets in
${\rm GL}_2^+(\RR)^2$, parameterised by $Q=(z,w)\in\HH^2$,
$$
Z_Q=\{(g,h)\in {\rm GL}_2^+(\RR)^2: gz, hw\in F {\rm\ and\ }
j(gz)+j(hw)=1\}.
$$
This family is definable in the o-minimal structure $\RR_{\rm an\ exp}$;
see e.g. [\PILAOAO].

Suppose that we have a solution $(x,y)$ to $(*)$ with large prime
$L=\max\{N,M\}$. Then we have $(u,v)$ with 
$\Phi_N(x,u)=0, \Phi_M(y,v)=0$.
Now $x,y$ cannot both be special, as $u+v=1$ contains no
special points (by K\"uhne [\KUHNE]). Let us assume for 
now neither is special.
So $\Phi_N(x,u)=0$, $\Phi_M(y,v)=0$,
and by the results of Najman [\NAJMAN] we have 
$$
[\QQ(u,v):\QQ]\ge cL^\delta.
$$

Take $z,w\in F$ with  $j(z)=x, j(w)=y$ and put $Q=(z,w)$.
Thus we have at least that many conjugate points $(u',v')$ over $\QQ$,
and each of these gives a solution of the system $(*)$ with the
same $(x,y)$. Each such $(u', v')$ gives rise to a rational point
on $Z_Q$, and (by results in [\HPA]) of height bounded by $CL^\eta$. 

By the Counting Theorem [\PILAWILKIE], which is uniform 
over the family, if
$L$ is sufficiently large then $Z_Q$ contains some positive 
dimensional real algebraic curve. The corresponding points
$(gz, hw)\in \HH^2$ must be non-constant, as the algebraic curves
in $Z_Q$ must account for ``many'' distinct $(u', v')$. So we
get a real algebraic curve contained in 
$$
\{(z,w)\in \HH^2: j(z)+j(w)=1\}.
$$
But then we must have a complex algebraic curve contained in it,
which then must coincide with it. This gives an algebraic curve
in $\HH^2$ whose image under $j$ in $\CC^2$ is algebraic. 
Then by the ``Ax-Lindemann'' theorem [\PILAOAO, Theorem 1.6], 
the image curve $u+v=1$ must be a modular curve. But it isn't. 
Thus $L$ is bounded.

Suppose $x$ is special. There are only finitely many rational
special points, so one is in a finite union of Hecke orbits, 
and for these one has a suitable Galois lower bound
(by isogeny estimates of Masser [\MASSER], 
subsequently refined by Pellarin [\PELLARIN] and others). 
So a similar argument applies, likewise
if $y$ is special.\ \qed

\medbreak
\noindent
{\bf Proof of Theorem 1.3.\/} This is exactly the same as above, 
except we appeal to SGH for $x\in\QQ$ instead of the results of 
[\NAJMAN] for the Galois lower bounds.\ \qed

\bigbreak
\centerline{\bf 3. How plausible is SGH?}

\medskip

SGH is related to uniform bounds for torsion in elliptic 
curves over number fields (Mazur, Kamienny, Merel,..[\MEREL]) 
and Serre's Uniformity Conjecture (Bilu-Parent,...[\BILUPARENT]) 
and seems in line with expectations.
A point $(x,y)\in T_N$ parameterises an elliptic curve 
with a cyclic subgroup of order $N$ defined over $\QQ(x,y)$. 
According to the Strong Uniform Boundedness Theorem 
of Merel (see [\MEREL, \REBOLLEDO]), 
the size of the torsion subgroup of $K$-rational points of
an elliptic curve defined over a numberfield $K$ with $[K:\QQ]=d$
is bounded by some $B(d)$. The known bounds for $B(d)$ are
exponential in $d$ but it is
conjectured that $B(d)$ can be taken polynomial in $d$
(see [\REBOLLEDO], Remark 2). 
The corresponding
conjectures for cyclic subgroups of size $N$,
i.e. for cyclic isogenies, would imply SGH.
The results of Najman [\NAJMAN] support these expectations.

That the gonality (defined in the proof of 7.3 below)
of modular curves grows at least as a positive 
power of $N$ is certainly necessary for SGH to hold. 
Conversely, Frey [\FREY] has shown 
(using Faltings's Big Theorem, i.e. his proof of Mordell-Lang
[\FALTINGSBIG]) that a curve has infinitely many points 
defined over fields
of degree $d$ over a field of definition $K$, then the gonality
of $C/K$ is at most $2d$. Thus, the modular curve $\Phi_N(x,y)=0$
has only finitely
many points defined over fields of degree at most $cN$ over $\QQ$
for some positive $c$.

\bigbreak
\centerline{\bf 4. Generalisation to curves}

\medskip
\noindent

There is nothing special about the curve $u+v=1$ in 
Theorems 1.1 and 1.3, except that it is not weakly special.
Both theorems hold for the system $(*)_V$ in which a 
non-weakly-special curve 
$V\subset \CC^2$ replaces the curve $u+v=1$ in $(*)$
and indeed 1.3 for $V$ is unconditional if $V$ is 
not defined over $\overline{\QQ}$. We do not formulate the results
as still more general formulations are in \S6.

If $V$ is special, say defined by $\Phi_K(u,v)=0$, then one 
can have rational solutions to $(*)_V$
with $x=y$ and arbitrarily large $\max\{N,M\}$. For if $\Phi_N(x, u)=0,
\Phi_M(x,v)=0$ then $u, v$ are Hecke equivalent, and one need
only choose $N,M$ such that this Hecke equivalence is given
by $\Phi_K$. Further, any weakly special curve whose fixed 
coordinate is in the Hecke orbit of a rational number will admit
rational solutions with arbitrarily large $\max\{N,M\}$ coming from the
non-fixed coordinate. But if we require $\min\{N,M\}\ge L$ then
only special subvarieties admit such points for arbitrarily large $L$
(under SGH for $x\in\QQ$ or unconditionally with $\max\{N,M\}$ prime.

This suggests the following ``Fermat-Mordell'' statement, in which a 
{\it weakly special subvariety\/} of $\GG_{\rm m}^k$ is a coset
of an algebraic subtorus. It is a consequence of the
$abc$ Conjecture (see e.g. [\HEIGHTS, Ch. 12] and below).

%

\medbreak
\noindent
{\bf 4.1. Conjecture.\/} Let $V\subset \GG_{\rm m}^2$ be a curve that
is not a weakly special subvariety. There is  $n(V)$ 
such that there are no rational points $(x^n, y^m)\in V,
x,y\in \QQ, x,y\ne 0,\pm1$ with $n,m\ge n(V)$.

\medbreak

We do not discuss here which multiplicative weakly special
varieties 
contain infinitely many such points. Some do.

Note that 4.1 is formulated in a slightly weaker form than the analogy
with 1.3 would suggest (which would be $\max\{n,m\}\ge n(V)$), 
in order to avoid possible issues if one exponent is small. 
This safer form is also adopted in subsequent conjectures.
One could formulate still more general conjectures addressing
solutions in the image of $(\QQ^\times)^2$ under 
morphisms $(\CC^\times)^2\rightarrow(\CC^\times)^2$ 
of large degree, or even correspondences, but this appears
to require some care and we defer this for now.

This conjecture clearly follow from Faltings's Theorem 
([\FALTINGS]) if the genus $g(V)\ge 2$
(with $n(V)=1$). If $g(V)\le 1$ the relation on $(x^n, y^m)$ could
still be of genus one or less for some small $n,m$.
These conjectures might be approachable for $V$ an elliptic curve. 

One can go further and state the following ``Fermat-Mordell-Lang'' 
formulation. Though apparently quite strong, it is nevertheless a 
consequence of the $abc$ Conjecture.

%
%
%

\medskip
\noindent
{\bf 4.2. Conjecture.\/} Let $\Gamma$ be a finitely generated
subgroup of $\QQ^\times$. There are only finitely many points 
$(u,v)=(sx^n, ty^m)$ on $u+v=1$ with $x,y\in \QQ$, $s,t\in\Gamma$, 
and $n,m\ge 4$. 

\medskip
\noindent
{\bf 4.3. Proposition.\/} {\it The $abc$ Conjecture implies 
Conjecture 4.2.\/}

\medskip
\noindent
{\bf Proof.\/} Let $\Gamma$ be a finitely generated subgroup of 
$\QQ^\times$. Enlarging if necessary, we may assume that $\Gamma$
is generated by $-1$ and some finite set $p_1,\ldots, p_k$
of prime numbers, and we set $P=p_1\ldots p_k$.
Now suppose we have a solution $(u,v)$ to the
equation in 4.2 with $n\ge m\ge 4$. Let us write $x=A/B, y=C/D$
where $A,B,C,D\in \ZZ$ are non-zero, with $(A,B)=(C,D)=1$.
By incorporating any $p_i$ that occur as factors into $s$ or $t$,
we may assume that $A,B,C,D$ are relatively prime to $P$
(and positive).
Multiplying through by a common denominator for $s, t$ and $B^n$
we have
$$
SA^n+TC^m{B^n\over D^m}=UB^n
$$
where $S, T, U$ are integers in $\Gamma$. 
We may assume they are relatively prime.
Since $(D,CT)=1$
we have $D^m\vert B^n$. Multiplying through by $D^m$ however
we conclude that $B^n\vert D^m$, hence they are equal.
So we have
$$
SA^n+TC^m=UB^n.
$$
The largest term in absolute value is either $TC^m$ or one of the terms
involving an $n$th power. Changing signs if needed, let us assume 
first that our equation is as above, with all terms positive.
By the $abc$ Conjecture (see e.g. [\HEIGHTS, Ch. 12]) with 
$\epsilon=1/4$
and $K=K_\epsilon$
we have
$$
UB^n<K{\rm rad}\big(SA^nTC^mUB^n\big)^{1+\epsilon}\le 
K\big(P \Big({U\over S}\Big)^{1/n}\Big({U\over T}\Big)^{1/n}
B^3\big)^{5/4}.
$$
Since $n\ge 4$ we find
$$
U^{3/8}B^{1/4}\le KP^{5/4}.
$$
Then $U,B$ are bounded, whence $S, T, A, C$ are also bounded.
The other case, when $TC^m$ is largest, is similar.\ \qed

\medskip

Of course one can also formulate a generalisation of 4.2 for
with a general (non-weakly-special) curve in place of $u+v=1$.
Note that the modular analogues of these do hold under 
SGH for $x\in\QQ$ (or unconditionally for tuples of isogenies 
where the largest degree is prime). 
That is because the notion of ``generation'' in the modular setting
is rather weak: the analogous 
statement is to seek points $(u,v): u+v=1$ where each of $u,v$
is either in the union of finitely many Hecke orbits or is in the
Hecke orbit of a rational number under a modular 
correspondence of large (prime) degree.

\bigbreak
\centerline{\bf 5. Generalisation to higher-dimensional varieties}

\medskip

The proof of Theorems 1.1 and 1.3 generalise without difficulty 
to higher  dimensions,  under the assumption of 
SGH for $x\in\QQ$ in generalising 1.3.

\medbreak
\noindent
{\bf 5.1. Definition.\/}
A {\it special subvariety\/} of $Y(1)^k$ is an irreducible component
of the intersection of (any number of) subvarieties of the following form:
$x_i=c$ where $c$ is constant and special; $\Phi(x_k, x_\ell)=0$
where $\Phi$ is a modular polynomial.
For a 
{\it weakly special subvariety\/},
the constant coordinates need not be special.
See e.g. [\HPA, \HPB, \PILAOAO].

\medskip
\noindent
{\bf 5.2. Theorem.\/} {\it 
Let $V\subset Y(1)^k$. There exists $L(V)$ with the following property.
If $u=(u_1,\ldots, u_k)\in V$ with $\Phi_{N_i}(x_i, u_i)=0$, 
$x_i\in \QQ$, $i=1,\ldots, k$,
and $N=\max\{N_i\}$
is a prime with $N\ge L(V)$ 
then $u$ lies in a 
positive dimensional weakly special variety contained in $V$.\/}
\medskip

\noindent
{\bf 5.3. Theorem.\/} {\it Assume SGH for $x\in \QQ$.
Let $V\subset Y(1)^k$. There exists $L(V)$ with the following property.
If $u=(u_1,\ldots, u_k)\in V$ with $\Phi_{N_i}(x_i, u_i)=0$,
$x_i\in \QQ$, $i=1,\ldots, k$,
and $\max\{N_i\}\ge L(V)$ then $u$ lies in a 
positive dimensional weakly special variety contained in $V$.\/}

\medskip
\noindent
{\bf Proof of 5.2 and 5.3.\/}
Let $K\subset\CC$ be finitely generated field of definition of $V$.
We take a definable family of sets 
$$
Z_Q=\{(g_1,\ldots, g_k)\in {\rm GL}_2^+(\RR)^k: g_iz_i\in F, i=1,\ldots, k, 
(j(g_1z_1),\ldots, j(g_kz_k)\in V\}
$$
parameterised by points $Q=(z_1,\ldots, z_k)\in \HH^k$.
Then a point $u=(u_1,\ldots, u_k)$ with 
$\Phi_{N_i}(x_i, u_i), i=1,\ldots, k$ and large
$L=\max\{N_i\}$ 
has (by the results of [\NAJMAN] for 5.2
and by SGH for $x\in{\QQ}$ for 5.3)
``many'' conjugates over $K$, and gives rise
to a $Q$ for which $Z_Q$  has ``many'' rational points.
By the Counting Theorem, we get a real algebraic arc in $Z_Q$
containing ``many'' of these points, 
and from it a real algebraic arc in 
$\{(z_1,\ldots, z_k): (j(z_1),\ldots, j(z_k))\in V\}$, hence a complex
algebraic curve contained there which, by the Ax-Lindemann
theorem  for the modular function [\PILAOAO],
is contained in a positive dimensional 
weakly special subvariety contained there, and it must be defined
over $\overline{K}$, as all coordinates of $u$ and its conjugates are.
The conjugates of this weakly special subvariety (over $K$) contain 
all the conjugates of $u$.
\ \qed

\medbreak

If one looks for points with large $\min\{N_i\}$, then 
(by an inductive argument) only special subvarieties
survive: under the same hypotheses and assumptions, there
is $L'(V)$ such that every point in $V$ of this form with 
$\min\{N_i\}\ge L'$ (and all $N_i$ is prime for the unconditional 
version) lies in a special
subvariety contained in $V$.

By analogy, one can formulate a conjectural generalisation of FLT 
in the setting of subvarieties of multiplicative groups.
As observed above, some 
weakly special subvarieties of $\GG_{\rm m}^k$ 
do have rational points which are arbitrarily large powers. 



\medskip
\noindent
{\bf 5.4. Conjecture.\/} Let $V\subset \GG_{\rm m}^k$. 
There is a positive integer $n(V)$ such that if
$P=(x_1^{n_1},\ldots, x_k^{n_k})\in V(\QQ)$, with all 
$x_i\in \QQ^\times, x_i\ne \pm1$ and
$n_i\ge n(V)$,
then $P$ lies in a positive dimensional 
weakly special subvariety of $\GG_{\rm m}^k$ contained in $V$.

\medskip

The  General Lang Conjecture ([\HEIGHTS, 14.3.7])
implies that all but finitely many
such points lie in the {\it special set\/} of $V$.
In the next section we will see that we can prove 5.4,
under a mild extra assumption, for the {\it inverse\/} Fermat equation.

Let SGH${}_d$ denote the special case of SGH in which 
$[\QQ(x):\QQ]\le d$. Under the assumption of SGH${}_d$,
the proofs of 1.3 and 5.3 go through if $x,y$
are restricted to be of degree at most $d$ over $\QQ$. One could 
then formulate all the conjectures above in this stronger form, 
with the hypothesis on the exponents now depending 
on $V, \Gamma, d$. 
The following conjecture is the most ambitious statement taking up
all these variants.

\medskip
\noindent
{\bf 5.5. Conjecture.\/} Let $V\subset \GG_{\rm m}^k$ be a 
subvariety defined
over $\CC$, let $\Gamma$ be a finite rank subgroup of $\CC^\times$, 
and let $d\ge 1$. There exists a constant
$n(V,\Gamma, d)$ with the following property. 
Suppose $P=(u_1,\ldots, u_k)\in V$ is a point such that, 
for $i=1,\ldots, k$, we have $u_i=s_ix_i^{n_i}$ with 
$s_i\in\Gamma$, $x_i$ not a root of unity, 
$[\QQ(x_i):\QQ]\le d$ and $n_i\ge n(V, \Gamma, d)$ then
$P$ lies in a positive-dimensional weakly special variety contained
in $V$.

\medbreak

It seems interesting to investigate whether Vojta's conjectures
(see e.g. [\HEIGHTS, Ch. 14]), which do imply Mordell-Lang, 
imply the above.

Let us conclude this section with a somewhat different 
generalisation of 1.3, and a further conjecture in the
multiplicative setting. We enunciate a different
weakening of SGH.

\medskip
\noindent
{\bf 5.6. Weak Galois-Orbit Hypothesis (WGH).\/} Let $F$
be a number field. There exists constants $c=c(F), \delta=\delta(F)$
such that if $(x,y)\in F\times\overline{\QQ}$ are not special
and $\Phi_N(x,y)=0$ then $[F(y):F]\ge cN^\delta$.

\medskip
\noindent
{\bf 5.7. Theorem.\/} {\it Assume WGH.
Let $K$ be a finitely generated subfield of $\CC$.
Let $V\subset Y(1)^k$.
There exists an integer $L=L(K, V)$ with the following property.
If $u=(u_1,\ldots, u_k)\in V$ with $\Phi_{N_i}(x_i, u_i)=0$,
$x_i\in K$, $i=1,\ldots, k$ and $\max\{N_i\}\ge L$ then
$u$ lies in a positive dimensional weakly special
subvariety contained in $V$.\/}

\medskip
\noindent
{\bf Proof.\/} We may assume that $V$ is defined over $K$.
Let $F=K\cap \overline{\QQ}$. Then $F$ is finitely
generated, and hence is a number field. 
Suppose $N_i=\max\{N_j, j=1,\ldots, k\}$.
For $x_i\in F$ and large $N_i$ the conclusion follows
using WGH and the proof of 5.3. For non-algebraic $x_i$
(and then also $u_i$) we apply Lemma 7.3 below to conclude
that $[K(y):K]\ge cN_i^\delta$ for suitable $c, \delta$ depending
on $K$, and then follow the proof of 5.3.\ \qed

\medskip

It is then natural to conjecture the analogous statement
in the multiplicative setting.

\medskip
\noindent
{\bf 5.8. Conjecture.\/}
Let $K$ be a finitely generated subfield of $\CC$.
Let $V\subset \GG_{\rm m}^k$.
There exists an integer $n=n(K, V)$ such that if
$P=(x_1^{n_1},\ldots, x_k^{n_k})\in V$, with 
$x_i\in K^\times$ but not a root of unity, 
and $n_i\ge n(V)$, for all $i=1,\ldots, k$,
then $P$ lies in a positive dimensional 
weakly special subvariety of $\GG_{\rm m}^k$ contained in $V$.

\medskip

I do not know  whether this statement in the case of plane curves
follows from the $abc$ conjecture. In the special case of
$V\subset \GG_{\rm m}^2$ defined by $u+v=1$, it asserts
the impossibility of solving
this equation in $K^{\times n}$, where $K$ is a finitely
generated field over $\QQ$, for large $n$ (depending on $K$).

\bigbreak
\centerline{\bf 6. Other settings and inverse Fermat}

\medskip
\noindent

It is natural to consider analogues in the setting of 
abelian varieties. 
The most natural analogue of 1.3 for an elliptic curve $E$
in place of $\GG_{\rm m}$
is the following statement, which is a consequence of
Mordell-Lang (ML; Faltings's Big Theorem [\FALTINGSBIG]) 
for $E\times E$. A {\it weakly special subvariety\/}
of an abelian variety is a translate of an abelian subvariety.

\medskip
\noindent
{\bf 6.1. A consequence of ML.\/} Let $E$ be an elliptic curve 
(defined over  $\CC$), and let $C\subset E\times E$ be a 
curve which is not weakly special.
There exists $L=L(E,C)$ with the following property. If $X,Y, U,V\in E$
are points such that: $U=[n]X, V=[m]Y, (U,V)\in C, X,Y\in E(\QQ)$
then $n,m\le L$.

\medbreak

Since $E(\QQ)$ is finitely generated, $U, V$ are in a finitely 
generated subgroup of $E\times E$. The statement then follows from
ML for $E\times E$.

One can consider a variant formulation.
Let $E$ be an elliptic curve in the form $y^2=x^3+ax+b$. 
Multiplication by $n$ on $E$ induces
an operation on $\CC$ as follows: $[n]x=z$ if $[n](x,y)=(z,w)$ on $E$.
There is a corresponding notion of ``weakly special'' variety in $\CC^2$,
comprising vertical and horizontal lines and the curves where
$[n]x=[m]y$ identically for some $n,m$.

\medskip
\noindent
{\bf 6.2. Conjecture.\/} Let $E$ as above and $V\subset \CC^2$
not ``weakly special''.
There exists $L=L(E,V)$  with the following property. If
$x,y\in \QQ$ and $([n]x, [m]y)\in V$ then $\max\{n,m\}\le L$.

\medskip

This statement is presumably not a consequence of ZP 
(as the points with rational $x$ are not finitely generated).

We now consider the analogue of conjecture 5.4 for the 
inverse Fermat equation:
 after all, $\Phi_n(x,u)=0$ is likewise the analogue
of $x=u^n$. On the inverse Fermat equation itself
see e.g. Lenstra [\LENSTRA].

Recall that, if $K$ is a field, $c\in K$, the polynomial $x^n-c$
is reducible over $K$ iff $c\in K^p$ for some prime number $p|n$,
or $c\in -4K^4$ and $4|n$ 
(see e.g. Lang [\LANGALGEBRA, VI, 9.1]).
The first condition is a natural minimality
for $u$ with $u^n=c\in K$: it guarantees that
$n$ is the {\it order\/}
of $u$ over $K$ in that no smaller power of $u$ lies in $K$.
The second condition
reflects the example $x^4+4=(x^2-2x+2)(x^2+2x+2)$.
Under the first condition only, one can get a lower bound on 
$[K(u):K]$, when
$K=\QQ$, from results of Risman [\RISMAN].

\medskip
\noindent
{\bf 6.3. Lemma.\/} {\it Let $\theta$ have order $n$ over $\QQ$.
Then $[\QQ(\theta):\QQ] \gg_\epsilon n^{1/2-\epsilon}$ for any 
$\epsilon>0$.\/}

\medskip
\noindent
{\bf Proof.\/} Write $h=[\QQ(\theta):\QQ]$. By [\RISMAN, Cor. 2],
we have $n=t\ell$ where $\ell$ divides $h$ and $\phi(t)$ divides $h$
(and $t$ is square-free). Either $t$ or $\ell$ must exceed $\sqrt{n}$.
\ \qed

\medskip
\noindent
{\bf 6.4. Theorem.\/} {\it Let $V\subset \GG_{\rm m}^k$. There is
a positive integer $n(V)$ with the following property.
Suppose $P=(u_1,\ldots, u_k)\in V$ and, for each $i=1,\ldots, k$,
$u_i^{n_i}=x_i$ where $x_i\in\QQ^\times$, and
$x_i\notin\QQ^p$ for any $p\vert n_i$. Suppose 
$\max\{n_1,\ldots, n_k\}\ge n$.
Then  $P$ lies in a positive-dimensional weakly special 
variety contained in $V$.\/}

\medskip
\noindent
{\bf Proof.\/} Under our assumptions, by Lemma 6.3,
the point $(u_1,\ldots, u_k)$ has degree at least 
$c\max\{n_1,\ldots, n_k\}^\delta$ over $\QQ$
for some absolute $c, \delta$, and hence will
have large degree over some fixed finitely generated field of definition
of $V$. Let $F=\RR\times[0, 2\pi]i$, a fundamental domain for the 
action of $2\pi i\ZZ$ on $\CC$ by translation. 
The restriction $\exp: F\rightarrow \CC^\times$ of the 
exponential function is definable in $\RR_{\rm an\ exp}$.
We take the definable family of sets
$$
Z_Q=\{(r_1,\ldots, r_k)\in \RR^k: z_j=2\pi i r_j\in F, j=1,\ldots, k,{\rm and}
$$
$$
\big(\exp(z_1+2\pi i r_1),\ldots, \exp(z_k+2\pi ir_k)\big)\in V\}.
$$
parameterised by points $Q=(z_1,\ldots, z_k)\in \CC^n$.
The rest of the proof is the same as the proof of 5.2 and 5.3,
using the Ax-Lindemann theorem for $\exp$ (a special
case of Ax-Schanuel [\AX]).\ \qed

\bigbreak
\centerline{\bf 7. Proof of Theorem 1.4}
\medskip

In this section we prove Theorem 1.4. The key point is that
modular curves have large gonality, and this implies that
transcendental points $(x,y): \Phi_N(x,y)=0$ give rise to extensions
of large degree over an arbitrary (but fixed) finitely generated 
extension of $\QQ$.
%
We first give a statement of the Zilber-Pink conjecture (ZP)
for subvarieties of $Y(1)^k$. 
See [\HPB] for various alternative formulations.

\medbreak
\noindent
{\bf 7.1. Definition.\/}
Let $V\subset Y(1)^k$. A subvariety $A\subset V$ is called
{\it atypical\/} (for $V$ in $Y(1)^k$) if there is a special subvariety
$T\subset Y(1)^k$ such that $A\subset V\cap T$ and
$$
\dim A> \dim V+\dim T-k.
$$



\medskip
\noindent
{\bf 7.2. Zilber-Pink Conjecture for $Y(1)^k$.\/} Let $V\subset Y(1)^k$.
Then $V$ has only finitely many maximal atypical subvarieties.

\medskip
\noindent
{\bf 7.3. Lemma.\/} {\it Let $K$ be a finitely generated subfield 
of $\CC$. There exist positive constants
$c, \delta$ (depending on $K$) with the following
property.
Let $P=(x,y)\in \CC^2$ be a point with non-algebraic coordinates
such that $\Phi_N(x,y)=0$. Then\/}
$$
[K(x,y):K]\ge cN^\delta.
$$

\medskip
\noindent
{\bf Proof.\/} Let us write $K=L(\kappa)$ where $L$ is a pure
transcendental extension of $\QQ$ and $[K:L]$ is a finite
algebraic extension. Do this minimising $[K:L]$ say.
Write $L=\QQ(t_1,\ldots, t_n)$ with the $t_i$
independent transcendental elements.

For a curve $C$ over a field $F$ with function-field $F(C)$ we write
$d_F(C)$ for its {\it gonality\/}: the minimum extension degree
$[F(C): F(t)]$ over $t\in F(C)$.

Let $P$ be such a point. We may assume that $x,y$ are algebraic 
over $K$. Let us choose $t_1,\ldots, t_m$ such that $x$ 
(and hence $y$)
are algebraic over $t_1,\ldots, t_m$ but not over 
$t_1,\ldots, t_{m-1}$. Let $M=\QQ(t_1,\ldots, t_{m-1})$ and
write $t=t_m$. The extension of fields $M(t,x,y)/ M(x,y)$ corresponds to
a dominant morphism of curves over $M$. Thus 
$$
d_M(M(t,x,y))\ge d_M(M(x,y))
$$
(see e.g. Poonen [\POONENGONALITY], where this fact is 
proved but described  as well-known).
Let $d_{\CC}(\Phi_N(x,y)=0)$ denote the $\CC$-gonality
of the modular curve. Then we have
$$
d_M(M(x,y))\ge d_{\CC}(\Phi_N(x,y)=0).
$$

Now $d_{\CC}(\Phi_N(x,y)=0)\ge c_0N$ for some positive 
constant $c_0$ (see [\ZOGRAF] and also [\ABRAMOVICH] 
where an explicit such bound is given).
Therefore
$$
[L(x,y):L]=[M(t,x,y): M(t)]  \ge c_0N,
$$
and so
$[K(x,y):K] \ge c_1N$
with $c_1=c_0/[K:L]$.
This proves the Lemma.\ \qed


\medbreak
\noindent
{\bf 7.4 Proof of Theorem 1.4.\/} For $A\subset Y(1)^k$, write $\langle A\rangle$ for the
smallest special subvariety of $Y(1)^k$ containing $A$.
If $\langle V\rangle\ne Y(1)^3$ then $V$ 
is atypical, and is then the unique maximal atypical subvariety. 
Conversely, if $V$ is atypical then it must be contained in a proper
special subvariety.
So we may assume that $V$ is not contained in any proper
special subvariety, and that atypical subvarieties of $V$ are points
which are contained in some special subvariety of codimension 2.

Suppose two coordinates, say $x,y$, are constant on $V$. 
They must be non-special and 
not in the same Hecke orbit. So a point
$(x,y,z)$ satisfying two special relations must be either a special
point $z$ that is in the Hecke orbit of either $x$ or $y$
(but then $x$ or $y$ would be special), or a $z$ which is
in the Hecke orbit of both $x$ and $y$ (but then $x$ and $y$ would be 
in the same Hecke orbit). Both are impossible.

Suppose just one coordinate, say $x$, is constant. Then the
image $V_{yz}$ of $V$ under projection to the $y,z$-plane
is a non-special curve, and we seek points which are either special
or in the Hecke orbit of $x$. Finiteness follows by
``Modular Mordell-Lang'' [\HPA, \PILAEMS].

So we may assume that no coordinate is constant on $V$.
We are looking for points $P=(x,y,z)$ satisfying two
special relations.
%
%
Let $P$ be such a point. 
It has one of the following forms: it is defined by two coordinates
being special; or by one coordinate being special and a modular
relation on the other two coordintates; or by modular relations
between two distinct pairs of coordinates.

Now if two coordinates are special, then we get a special point
on the image of $V$ under projection to those coordinates.
This image is not special (since $\langle V\rangle =Y(1)^3$),
and so for each choice of pair of coordinates there are only
finitely many such points.

If $P$ is a point of the second type, we distinguish two subcases.
In the first subcase, the two modular related points are algebraic.
Then $P$ is an algebraic point of $V$ and in a finite set. In the second
subclass, the two modular-related coordinates are transcendental
over $\overline{\QQ}$. Such $P$ then has ``many'' conjugates
over $K$, by a combination of Lemma 7.3 and Landau-Siegel.
We conclude this case by o-minimality and point-counting,
much as we deal with the following final case.

The last case concerns points $P$ satisfying modular relations on two
distinct sets of coordinates. So all three coordinates of $P$ are
in the same Hecke orbit. By Lemma 7.3, $P$ 
has ``many'' conjugates over $K$,
and thus $V$ contains ``many'' points $P'$ which are intersections
with special subvarieties of the same complexity as the one
containing $P$.

Let $Z\subset \HH^3$ be the preimage of $V$ in $\HH^3$
intersected with $F^3$, where $F$ is the standard fundamental 
domain for the action of ${\rm SL}_2(\ZZ)$. Then $Z$
is definable. For $g,h\in {\rm GL}_2^+(\RR)$ we have the
Mobius subvariety $M_{g,h}\subset \HH^3$ defined by
$$
M_{g,h}=\{(u,gu, hgu)\in \HH^3: u\in \HH\}.
$$
We consider the following definable subset of
${\rm GL}_2^+(\RR)^2$:
$$
W=\{(g,h): M_{g,h}\cap Z \neq \emptyset\}.
$$

Each conjugate $P'$ of $P$ over $K$ gives rise to a rational point
$(g,h)\in W$ whose height is $\le c\langle P\rangle ^C$, and
we get $\ge c'\langle P\rangle ^{C'}$ such points.
By the Counting Theorem, $W$ contains positive-dimensional semi-algebraic sets, and the intersection points of the corresponding Mobius subvarieties with $Z$ must move, by the same argument used in [\HPA],
in order to account for the ``many'' distinct pre-images of the $P'$.

Complexifying the real parameter of the moving family of Mobius subvarieties we get a complex surface in $\HH^3$ which intersects
$Z$ in a set of at least one real dimension, and hence in a set of
one complex dimension, and so contains the premiere of $V$.
By Ax-Schanuel ([\PTAXJ], though in this case in fact just the 
special case
``Ax-Logarithms'' established in [\HPA]), $V$ is contained in a proper
weakly special subvariety of $\CC^3$.

But this is a contradiction, as $V$ is not contained in a proper
special subvariety (by hypothesis), and no coordinate is constant
on $V$ (as we reduced to this case). \ \qed

\medbreak
\noindent
{\bf 7.5. Proposition.\/}
{\it Let $V\subset Y(1)^4$ be a curve which is not contained in any proper special subvariety and assume that no image 
of $V$ under a coordinate projection to $Y(1)^3$ is defined over
$\overline{\QQ}$. Then there are only finitely many points
$(w,x,y,z)\in V$ such that, for some $N,M$,
$\Phi_N(w,x)=0, \Phi_M(y,z)=0$.\/}

\medbreak
\noindent
{\bf Proof.\/} Suppose two coordinates are constant on $V$. 
Say $w$ is one of them.
If $x$ is also constant we cannot have
$\Phi_N(w,x)=0$, for then $V$ would be contained in a proper
special; and for other $x$ there are no points of the required form.
If, say, $z$ is also constant then $x,y$ are non-constant (by above)
and $V$ projects to a curve $V_{xy}$ in the $xy$-plane.
We are looking for points in $V_{xy}$ whose $x,y$ coordinates are
in the Hecke orbits of $w,z$, respectively, and finiteness
follows by Modular Mordell-Lang as above.

Suppose just one coordinate, say $w$ is constant. So $y,z$ are
non-constant and satisfy some algebraic relation. If this relation is
not defined over $\overline{\QQ}$ then, with finitely many exceptions,
the sought points have $y,z$ non-algebraic. Let $K$ be a finitely 
generated field of definition of $V$. Then $[K(x):K]\ge cN^\delta$ 
for some $c, \delta>0$ by isogeny estimates, and 
$[K(y,z):K] > cM^\delta$ for some $c,\delta>0$ by gonality,
and an argument using o-minimality, point-counting and
Modular Ax-Lindemann concludes as above.

So we can suppose that no coordinates are constant on $V$,
and so every pair of coordinates satisfy some algebraic relation.
Suppose neither of the relations $R(w,x)=0, S(y,z)=0$ are
defined over $\overline{\QQ}$. Then, with finitely many exceptions,
each pair $(w,x), (y,z)$ consists of transcendental points.
These have large degree over $K$ in relation to the complexity
$\max\{N,M\}$, and we conclude as above. 

If on the other hand both these pairs of relations are over 
$\overline{\QQ}$ then $w,x,y,z$ are all algebraic, and there are 
only finitely many points when even three of the coordinates are,
under our hypotheses.

We are reduced to the case that $R(w,x)=0$, say, is defined
over $\overline{\QQ}$, but $S(y,z)=0$ is not. 
Consider the curve image $V_{xyz}$ under projection to the
$xyz$ coordinates. We are looking for points where $x$ is algebraic,
and $y,z$ have a modular relation. If $V_{xyz}$ is not contained
in a proper subvariety of $Y(1)^3$ defined over $\overline{\QQ}$
then the finiteness of such $(x,y,z)$ is a trivial consequence
of the main theorem of [\CGMM]. 

So we may assume that
$V_{xyz}$ is contained in a proper subvariety $W$ defined over 
$\overline{\QQ}$, defined say by $P(x,y,z)=0$. We observe that
there can be only finitely many $x$ for which the relation $P(x,y,z)$
on $y,z$ is divisible by modular relation. 
For other $x$, if $x$ is algebraic
and $\Phi_M(y,z)=0$ then $y,z$ must also be algebraic, 
and there are only finitely many such points.\ \qed

\medbreak
\noindent
{\bf 7.6. Theorem.\/} {\it Let $V\subset Y(1)^k$ be a curve such
that no image of it under projection to three coordinates is defined
over $\overline{\QQ}$. Then ZP holds for $V$.\/}

\medbreak
\noindent
{\bf Proof.\/}
As above, we may assume that $V$ is not contained in any proper
special subvariety of $Y(1)^k$. We consider atypical points,
and these involve either 2 coordinates (for two points being special),
or 3 coordinates, or 4 coordinates (the case of modular
correspondences between disjoint pairs of coordinates). 
In each case, finiteness
is covered by either 7.4 or 7.5.\ \qed

\bigbreak
\centerline{\bf 8. SGH and GO1}

\medskip

In this section we show that SGH in fact implies the statement
formulated as GO1 in [\HPB], of which it is a special case.

We define the {\it complexity\/} of a special subvariety as follows.
If $x\in \CC$ is special,
we denote by $D(x)$ the discriminant of the corresponding quadratic 
order (i.e. the endomorphism ring of the elliptic curve $E$
with $j$-invariant $x$). Alternatively, $D(x)$ is the discriminant
$b^2-4ac$ where $az^2+bz+c=0$ is the minimal polynomial
of some pre-image $z=j^{-1}(x)$ of $x$ over $\ZZ$.

\medskip
\noindent
{\bf 8.1. Definition.\/} The {\it complexity\/} of a special subvariety
$T\subset Y(1)^k$ is 
$$
\Delta(T)=\max\{ D(x_i), N(x_h, x_\ell)\}
$$
where $D(x_i)$ ranges over all constant coordinates, and
$N(x_h, x_\ell)=N$ if $x_h, x_\ell$ are non-constant coordinates 
which are related by a modular polynomial $\Phi_N$, 
and we range over all such related pairs.

\medskip
\noindent
{\bf 8.2. Formulation GO1 ([\HPB]).\/} 
Let $V\subset Y(1)^k$ be defined over a field $K$ which
is finitely generated over $\QQ$. There are positive constants
$c, \eta$ with the following property. If $P\in V$ defined over
a field extension of $K$ then
$
[K(P):K]\ge c\Delta(\langle P\rangle)^\eta.
$

\medskip
\noindent
{\bf 8.3. Proposition.\/} {\it SGH implies GO1 for the subvarieties
$Y(1)^k\subset Y(1)^k, k=1,2,\ldots$ as a subvarieties defined over 
$\QQ$.\/}

\medskip
\noindent
{\bf Proof.\/} Suppose $x=(x_1,\ldots, x_k)\in Y(1)^k$. 
Some $x_i$ may  be special, and some pairs of coordinates 
may be  related by modular polynomials. For the special 
$x_i$ we have Landau-Siegel. 
Suppose $x_{i_1},\ldots, x_{i_k}$ are all 
in the same Hecke orbit. The
complexity of $\langle (x_{i_1},\ldots, x_{i_k})\rangle$
is then the maximum $N$ of the $N_{ab}$
such that $\Phi_{N_{ab}}(x_{i_a}, x_{i_b})=0$, and by SGH
we have $[\QQ(x_1,\ldots, x_n):\QQ]\ge cN^\delta$.\ \qed

\medskip
\noindent
{\bf 8.4. Proposition.\/} {\it GO1 for 
$Y(1)^k\subset Y(1)^k, k=1,2,\ldots$
implies GO1 in general.\/}

\medskip
\noindent
{\bf Proof.\/}
Assume the truth of GO1 for 
$Y(1)^n\subset Y(1)^n, n=1,2,\ldots$ and let $V\subset Y(1)^n$
defined over a field $K$ finitely generated over $\QQ$.
Let us write $K=L(\kappa)$ where $L$ is purely
transcendental over $\QQ$ and $[K:L]$ is algebraic.
Let $P=(x_1,\ldots, x_n)\in V$. We may suppose all coordinates
are algebraic over $K$.

Some coordinates of $P$ may be special, and some related by
modular polynomials. If $x_i$ is special, then it is algebraic
and its degree over $\QQ$ is bounded below by 
$c\Delta(x)^\delta$ be Landau-Siegel. 
If $\Phi_N(x_i, x_j)$ then we distinguish two cases. If one (and hence
both) $x_i, x_j$ are algebraic, the required lower bound follows
from SGH. If they are not algebraic, then 
$[\QQ(x_i, x_j): \QQ(x_i)]=\deg \Phi_N$, and the required degree
bound follows via the gonality argument in 
the proof of 7.2.\ \qed

\medskip

Note that GO1, is stronger than 
the conjectured ``LGO'' used in [\HPB] to give
a conditional proof of the Zilber-Pink conjecture for $Y(1)^k$
(a second condition in [\HPB], a suitable ``Ax-Schanuel'' 
statement for the modular function,
has subsequently been affirmed in [\PTAXJ]).
Thus SGH implies the full Zilber-Pink conjecture for $Y(1)^k$.

\bigbreak
\noindent
{\bf Acknowledgements.\/} 
I am very grateful to Peter Sarnak for several 
comments and suggestions.
The work was partially supported by 
the EPSRC, grant  EP/J019232/1.

\bigbreak

\noindent
{\bf References\/}

\medskip

\item{\ABRAMOVICH.\/} D. Abramovich,
A linear lower bound for the gonality of modular curves,
{\it IMRN\/} {\bf 1996}, No. 20.

\item{\ANDREBOOK.} Y. Andr\'e, {\it $G$-functions and geometry,\/} 
Aspects of Mathematics E13, Vieweg, Braunschweig, 1989.

\item{\ANDRE.} Y. Andr\'e, Finitude des couples d'invariants 
modulaires singuliers sur une courbe alg\'ebrique plane non modulaire, 
{\it Crelle\/} {\bf 505} (1998), 203--208.

\item{\AX.} J. Ax, On Schanuel's
conjectures, {\it Annals\/} {\bf 93\/} (1971), 252--268.

\item{\BILUPARENT.} Y. Bilu and P. Parent,
Serre's uniformity problem in the split Cartan case,
{\it Annals\/} {\bf 173} (2011), 569--584.

\item{\HEIGHTS.}
E. Bombieri and W. Gubler, {\it Heights in Diophantine Geometry},
New Mathematical Monographs {\bf 8},
Cambridge University Press, 2006.

\item{\CGMM.} Z. Chatzidakis, D. Ghioca, D. Masser, and
G. Maurin, Unlikely, likely and impossible intersections
without algebraic groups, 
{\it Rend. Link Mat. Apple.\/} {\bf 24} (2013), 485--501.

\item{\FALTINGS.} G. Faltings,
Endlichkeitss\"atze f\"ur abelsche Variet\"aten \"uber
Zahlk\"orperen,
{\it Inventiones\/} {\bf 73} (1983), 349--366.

\item{\FALTINGSBIG.} G. Faltings,
Diophantine approximation on abelien varieties,
{\it Annals\/} {\bf 133} (1991), 549--576.

\item{\FREY.} G. Frey,
Curves with infinitely many points of fixed degree,
{\it Israel J. Math.\/} {\bf 85} (1994), 79--83.

\item{\HPA.} P. Habegger  and  J. Pila, Some
unlikely intersections beyond  Andr\'e--Oort, 
{\it Compositio\/} {\bf 148} (2012), 1--27.

\item{\HPB.} P. Habegger and J. Pila, 
O-minimality and certain atypical intersections, arXiv,
and {\it Annales de l'ENS}, to appear.

\item{\KLINGLERYAFAEV.} B. Klingler and A. Yafaev, 
The Andr\'e-Oort conjecture, {\it Annals\/} {\bf 180} (2014), 867--925.

\item{\KUHNE.} L. K\"uhne, An effective result of Andr\'e-Oort type,
{\it Annals\/} {\bf 176} (2012), 651--671.

\item{\LANG.} S. Lang, Integral points on curves,
{\it Publ. Math. IHES\/} {\bf 6} (1960), 319--335.

\item{\LANGALGEBRA.\/} S. Lang, {\it Algebra\/}, Revised Third Edition,
GTM 211, Springer, New York, 2002.

\item{\LAURENT.} M. Laurent, \'Equations diophantiennes 
exponentielles, {\it Invent. Math.\/} {\bf 78} (1984), 299--327.

\item{\LENSTRA.\/} H. W. Lenstra, On the inverse Fermat equation,
{\it Discrete Mathematics\/} {\bf 106/107} (1992), 329--331.

\item{\LIARDET.} P. Liardet, Sur une conjecture de Serge Lang,
{\it Asterisque} {\bf 24-25} (1975), 187--210.

\item{\MANN.} H. B. Mann, On linear relations between roots of unity, 
{\it Mathematika\/} {\bf 12} (1965), 107--117.

\item{\MASSER.} D. Masser,
Small values of the quadratic part of the N\'eron-Tate height
on an abelian variety,
{\it Compositio\/} {\bf 53} (1984), 153--170.

\item{\MEREL.} L. Merel, Bornes pour la torsion des 
courbes elliptiques sur les corps de nombres,
{\it Inventiones\/} {\bf 124\/} (1996), 437--449.

\item{\NAJMAN.} F. Najman, Isogenies of non-CM elliptic curves with
rational $j$-invariants over number fields,
arXiv:1506.03127.

\item{\OORT.} F. Oort, Canonical lifts and dense sets of CM points, {\it Arithmetic Geometry, Cortona, 1994,\/} 228--234, F. Catanese, editor,
Symposia. Math., XXXVII, CUP, 1997.

\item{\PELLARIN.}
F. Pellarin, Sur une majoration explicite pour un degr{\'e}
d'isog{\'e}nie liant deux courbes elliptiques, 
{\it Acta Arith.\/} {\bf 100} (2001), 203--243.

\item{\PILAOAO.} J. Pila, 
O-minimality and the Andr\'e-Oort conjecture 
for $\CC^n$, {\it Annals\/} {\bf 173} (2011), 1779--1840.

\item{\PILAEMS.} J. Pila, 
Special point problems with elliptic modular surfaces, 
{\it Mathematika\/} {\bf 60} (2014), 1--31.

\item{\PILAICM.} J. Pila,
O-minimality and Diophantine geometry,
{\it Proceedings of the ICM, Seoul, 2014\/},
S. Y. Jang, Y. R. Kim, D.-W. Lee, and I. Yie, eds.,
Vol. I, 547--572,
Kyang Moon SA, Seoul, 2014.

\item{\PTAG.} J. Pila and J. Tsimerman,
Ax-Lindemann for ${\cal A}_g$, 
{\it Annals\/} {\bf 179} (2014), 659--681.

\item{\PTAXJ.} J. Pila and J. Tsimerman,
Ax-Schanuel for the $j$-function, arXiv, submitted.

\item{\PILAWILKIE.} J. Pila and A. J. Wilkie, 
The rational points of a definable set, 
{\it DMJ\/} {\bf 133} (2006), 591--616.

\item{\PILAZANNIER.} 
J. Pila and U. Zannier, Rational points in periodic analytic sets and the 
Manin-Mumford conjecture, {\it Rend. Lincei Mat. Appl.} {\bf 19} (2008), 
149--162.

\item{\PINK.} R. Pink, A common generalization of the 
conjectures of Andr\'e-Oort, Manin-Mumford, and Mordell-Lang, 
manuscript dated 17 April 2005,  available from 
{\tt http://www.math.ethz.ch/$\sim$pink/.}

\item{\POONENGONALITY.} B. Poonen,
Gonality of modular curves in characteristic $p$,
{\it Math. Res. Lett.} {\bf 14} (2007), 691--701.

\item {\RAYNAUDA.} M. Raynaud, 
Courbes sur une vari\'et\'e ab\'elienne
et points de torsion, {\it Invent. Math.} { {\bf 71} (1983), no. 1, 207--233.

\item{\RAYNAUDB.}  M. Raynaud, 
Sous-vari\'et\'es d'une vari\'et\'e ab\'elienne et
points de torsion, in {\it Arithmetic and Geometry, Volume I, \/} 
pp 327--352, 
Progr. Math. {\bf 35}, Birkhauser, Boston MA, 1983.

\item{\REBOLLEDO.} M. Rebolledo,
Merel's theorem on the boundedness of the torsion of
elliptic curves, {\it Clay Mathematics Proceedings\/} {\bf 8} (2000), 
71--82.

\item{\RISMAN.} L. J. Risman,
On the order and degree of solutions of pure equations,
{\it Proc. Amer. Math. Soc.\/} {\bf 55} (1976), 261--266.

\item{\SARNAKADAMS.} P. Sarnak and S. Adams, 
Betti numbers of congruence subgroups (with an appendix by 
Z. Rudnick), {\it Israel J. Math.} {\bf 88} (1994), 31--72.

\item{\SIEGEL.\/} C.-L. Siegel, 
\"Uber einige Anwendungen diophantischer
Approximationen, {\it Abh. Preuss. Akad. Wissen. Phys. Math.
Kl.\/} (1929), 41--69. Also, {\it Gesammelte Abhandlungen\/}, Vol. I.
English translation by C. Fuchs in 
{\it On Some Applications of Diophantine Approximations,\/}
U. Zannier, ed., Edizioni della Normale, 2014.

\item{\TSIMERMAN.} J. Tsimerman,
A proof of the Andr\'e-Oort conjecture for ${\cal A}_g$,
arXiv:1506.01466.

\item{\ULLMOYAFAEV.} E. Ullmo and A. Yafaev, Galois orbits 
and equidistribtuion of special subvarieties: towards the Andr\'e-Oort
conjecture, {\it Annals} {\bf 180} (2014), 823--865.

\item{\WILES.\/} A. Wiles,
Modular elliptic curves and Fermat's Last Theorem,
{\it Annals\/} {\bf 141} (1995), 553--572.

\item{\WILKIE.\/} A. J. Wilkie,
Model completeness results for
expansions of the ordered field of real numbers by restricted
Pfaffian functions and the exponential function, {\it J.
Amer. M. Soc.\/} {\bf 9} (1996), 1051--1094.

\item{\ZANNIER.} U. Zannier, 
{\it Some Problems of Unlikely Intersections
in Arithmetic and Geometry,\/} with appendices by D. Masser,
{\it Annals of Mathematics Studies\/} {\bf 181}, 
Princeton University Press, 2012.

\item{\ZILBER.} B. Zilber, Exponential sums equations and the Schanuel 
conjecture, {\it J. London Math. Soc. (2)\/} {\bf 65} (2002), 27--44.

\item{\ZOGRAF.} P. G. Zograf, Small eigenvalues of automorphic 
Laplacians in spaces of parabolic forms, {\it Zap. Nauchn. Sem. 
Leningrad Otdel. Mat. Inst. Steklov (LOMI)\/} {\bf 134}
(1984), 157--168. In Russian; English translation in 
{\it J. Soviet Math.\/} {\bf 36} (1987), 106--114.

\vfill

\noindent
Mathematical Institute

\noindent
University of Oxford

\noindent
Oxford OX2 6GG 

\noindent
UK

\bye